\documentclass[aos,MSNbibl,nameyear,seceqn,dvips]{arximspdf}

\doi{10.1214/11-AOS888}
\volume{39}
\issue{4}
\pubyear{2011}
\firstpage{1852}
\lastpage{1877}

\makeatletter
\newtheorem{theo}{Theorem}[section]
\newtheorem{lemma}[theo]{Lemma}
\newtheorem{fact}{Fact}
\makeatother

\begin{document}
\begin{frontmatter}

\title{Optimal model selection for density estimation of stationary
data under various mixing conditions}
\runtitle{Optimal model selection under various mixing conditions}

\begin{aug}
\author[A]{\fnms{Matthieu} \snm{Lerasle}\corref{}\ead[label=e1]{lerasle@gmail.com}\thanksref{T1}}
\runauthor{M. Lerasle}
\thankstext{T1}{Supported in part by FAPESP Grant 2009/09494-0.}
\affiliation{IMT (UMR 5219), INSA Toulouse}
\address[A]{IMT (UMR 5219), INSA Toulouse\\
135 Avenue de Rangueil\\
31077 Toulouse Cedex 4\\
France\\
\printead{e1}} 
\end{aug}

\received{\smonth{2} \syear{2011}}

%
\begin{abstract}
We propose a block-resampling penalization method for marginal density
estimation with nonnecessary independent observations.\break When the data
are $\beta$ or $\tau$-mixing, the selected estimator satisfies oracle
inequalities with leading constant asymptotically equal to $1$.

We also prove in this setting the slope heuristic, which is a
data-driven method to optimize the leading constant in the penalty.
\end{abstract}

%
\begin{keyword}[class=AMS]
\kwd[Primary ]{62G07}
\kwd{62G09}
\kwd[; secondary ]{62M99}.
\end{keyword}
\begin{keyword}
\kwd{Density estimation}
\kwd{optimal model selection}
\kwd{resampling methods}
\kwd{slope heuristic}
\kwd{weak dependence}.
\end{keyword}
\vspace*{6pt}
\end{frontmatter}
%

\section{Introduction}
Model selection by penalization of an empirical loss is a~general
method that includes several famous procedures as cross-validation
[\citet{Ru82}] or hard thresholding [\citet{DJKP}] as shown by \citet{BBM99}.
The difficulty is to calibrate the penalty so that the
selected estimator satisfies an oracle inequality. A good penalty has
the shape of an ideal one [see definition (\ref{defPenid})] and
depends in general on a leading constant that should be chosen
sufficiently large.

Resampling penalties provide a shape for the penalty term in a general
statistical learning framework; see \citet{Ar08}. The resulting
estimator satisfies sharp oracle inequalities in non-Gaussian
heteroscedastic regression among histograms [\citet{Ar08}] and in
density estimation among more general collections of models [\citet{Le09}].
The validity of these theorems relies on the independence of
the observations. In this paper, we study a~generalization of these
penalties, called block-resampling penalties and we prove that the
resulting estimator satisfies sharp oracle inequalities when the data
are only supposed to be $\beta$- or $\tau$-mixing [the coefficient
$\beta$ has been defined by \citet{RV59}, the coefficient $\tau$ by
\citet{DP05}; see Section \ref{Ch4ss23}].

We use a coupling method to extend the results for independent data. It
was introduced in \citet{BCV01} in a regression problem and used in
\citet{CM02} for density estimation with $\beta$-mixing observations. $\beta$
is a well known ``strong'' mixing coefficient. We refer to the books of
\citet{Do94} and \citet{Br07} for examples of $\beta$-mixing processes.
One of the most important is the following: a~stationary, irreducible,
aperiodic and positively recurrent Markov chain is $\beta$-mixing.
``Strong'' mixing coefficients cannot be used to study a lot of simple
processes. For example, the stationary solution of the equation
%
%
\begin{equation}\label{eqan84}
X_{n}=\tfrac12(X_{n-1}+\xi_n),
\end{equation}
where $(\xi_n)_{n\in\mathbb{Z}}$ are i.i.d. Bernoulli random variables
$\mathcal{B}(1/2)$ is not $\beta$-mixing [see \citet{An84}]. This is why
``weak'' mixing coefficients such as $\tau$ have been introduced. They are
easier to compute and allow us to cover more examples, as the process
(\ref{eqan84}) [we refer to the papers of \citet{DP05}, \citet{CDT08} or the
book of \citet{DDLLLP} for examples of weakly-mixing processes]. In
\citet{Le08}, we used a~coupling result of \citet{DP05} to extend the coupling
method to $\tau$-mixing data.

In all these previous papers, the dimension of the models was used as a
shape of the penalties. The leading constant was built with the mixing
coefficients and could not in general be computed from the data. When
it could, the theoretical upper bounds obtained are probably too
pessimistic to be used by the statistician. We use in this paper
block-resampling penalties as a shape of the penalty and a data-driven
leading constant. The first main result of the paper is that the
resulting estimator satisfies asymptotically optimal oracle inequalities.
We propose also to optimize the choice of the leading constant in
penalties using the slope algorithm. This procedure is based on the
slope heuristic, introduced in \citet{BM07} and proved in \citet{BM07}
for Gaussian regression, in \citet{AM08} for non-Gaussian
heteroscedastic regression over histograms and in \citet{Le09} for
density estimation. The second main result of this paper is a proof of
the slope heuristic for the marginal density estimation problem with
$\beta$- or $\tau$-mixing data.

Block-resampling penalties and the slope heuristic can be defined in
a~more general statistical learning framework, including the problems of
classification and regression [see \citet{Ar08}; \citet{AM08}]. Our results are
contributions to the theoretical understanding of these generic
methods. Up to our knowledge, they are the first ones obtained in a
mixing framework.

The paper is organized as follows. Section \ref{Ch4S2} introduces the
density estimation framework, the estimators, the penalties and the
main assumptions. Sections~\ref{Ch4S3} and~\ref{Ch4S4} give the main
results, respectively, for $\tau$- and $\beta$-mixing processes. Section~\ref{Ch4S5}
gives the proofs of the main results. Some other proofs are
available as Supplementary Material [\citet{le11}].
%
\section{Preliminaries}\label{Ch4S2}
%
\subsection{The density estimation framework}
We observe $n$ real valued, identically distributed random variables
$X_1,\ldots,X_n$, defined on a probability space $(\Omega,\mathcal
{A},\mathbb{P}
)$, with common law $P$. We assume that $P$ is absolutely continuous
with respect to the Lebesgue measure $\mu$ on $\mathbb{R}$ and we
want to
estimate the density $s$ of $P$ with respect to $\mu$. $L^2(\mu)$
denotes the Hilbert space of square integrable real valued functions
and $\|\cdot\|$ the associated $L^2$-norm. We assume that $s$ belongs to
$L^2(\mu)$. The risk of an estimator $\widehat{s}$ of $s$ is measured
with the
$L^2$-loss, that is $\|s-\widehat{s}\|^2$, which is random when
$\widehat{s}$ is.

Let $p$, $q$ be two integers and assume that $n=2pq$. For all
$i=0,\ldots,p-1$, let $I_i=(2iq+1,\ldots,(2i+1)q)$, $A_i=(X_l)_{l\in I_i}$.
For all functions $t$ in $L^1(P)$, for all reals $x_1,\ldots,x_q$, we define
\begin{eqnarray*}
L_qt(x_1,\ldots,x_q)&=&\frac1q\sum_{i=1}^qt(x_i),\qquad Pt=\int_\mathbb
{R}t(x)s(x)\,d\mu
(x),\\
 P_At&=&\frac1p\sum_{i=0}^{p-1}L_qt(A_i).
\end{eqnarray*}
Given a linear space $S_m$ of measurable, real valued functions, and an
orthonormal basis $(\psi_{\lambda})_{\lambda\in\Lambda_m}$ of
$S_m$, we define the projection
estimator $\widehat{s}_{A,m}$ of~$s$ onto $S_m$ by
\[
\widehat{s}_{A,m}=\sum_{\lambda\in\Lambda_m}(P_A\psi_{\lambda
})\psi_{\lambda}\in\arg\min_{t\in S_m}\{ \|t\| ^2-2P_At \}.
\]
Given a finite collection $(S_m)_{m \in\mathcal{M}_n}$ of such linear
spaces and a
penalty function $\operatorname{pen}\dvtx \mathcal{M}_n\rightarrow\mathbb
{R}_+$, the Penalized Projection
Estimator, hereafter PPE, is defined by
%
%
\begin{equation}\label{defPPE}
\qquad\tilde{s}_A=\widehat{s}_{A,\widehat{m}} \qquad\mbox{where } \widehat{m}\in\arg
\min_{m \in\mathcal{M}_n
}\{ \|\widehat{s}_{A,m}\|^2-2P_A\widehat{s}_{A,m}+\operatorname
{pen}(m) \}.
\end{equation}
We will say that the PPE satisfies an oracle inequality when one of the
two following inequalities holds.

There exist constants $ \kappa>0, \gamma>1$ and a positive sequence
$(K_n)_{n\in\mathbb{N}^*}$ bounded away from $0$ such that
%
%
\begin{equation}\label{defOItraj}
\mathbb{P}\Bigl(K_n\|s-\tilde{s}_A\|^2\leq\inf_{m \in\mathcal{M}_n}\|
s-\widehat{s}_{A,m}\|^2
\Bigr)\geq1-\frac{\kappa}{n^{\gamma}}.
\end{equation}
There exists a positive sequence $(K_n)_{n\in\mathbb{N}^*}$ bounded
away from
$0$ such that
%
%
\begin{equation}\label{defOImoyenne}
K_n\mathbb{E}(\|s-\tilde{s}_A\|^2)\leq\mathbb{E}\Bigl(\inf_{m \in
\mathcal{M}_n}\|s-\widehat{s}
_{A,m}\|^2\Bigr).
\end{equation}
The oracle inequality is said to be sharp when, moreover, the sequence~\mbox{$K_n\!\rightarrow\!1$}
when $n$ grows to infinity. Inequalities (\ref{defOItraj}) are usually preferred to (\ref{defOImoyenne}) since they
describe the typical behavior of the selected estimator and not only of
its expectation.

It is worth mentioning that we only use $\operatorname{Card}(\bigcup_{i=0}^{p-1}
I_i)=pq=n/2$ data to build the estimator $\tilde{s}_{A}$. The
consequences of this choice are discussed after Theorem \ref{theoRestau} and in Section \ref{Ch4ss43}.
%
\subsection{Block-resampling penalties}
We introduce block-resampling penalties as natural generalizations of
resampling penalties. The best estimator in the collection $(\widehat{s}
_{A,m})_{m \in\mathcal{M}_n}$ minimizes among $\mathcal{M}_n$ the
ideal criterion
\[
\|s-\widehat{s}_{A,m}\|^2-\|s\|^2=\|\widehat{s}_{A,m}\|
^2-2P_A\widehat{s}_{A,m}+\operatorname{pen}_{\mathrm{id}}
(m).
\]
In this decomposition, the ideal penalty $\operatorname{pen}_{\mathrm
{id}}(m)$ [\citet{Ar08}] is
equal to
%
\begin{equation}
\operatorname{pen}_{\mathrm{id}}(m)=2(P_A-P)(\widehat
{s}_{A,m}).\label{defPenid}
\end{equation}
To adapt the approach of \citet{Ar08} to a dependent setting, we replace
the resampling step by a resampling procedure on the blocks
$(A_i)_{i=0,\ldots,p-1}$. Let $(W_0,\ldots,W_{p-1})$ be a resampling scheme,
that is, a vector of positive random variables, independent of
$(X_i)_{i=1,\ldots,n}$ and exchangeable, which means that, for all
permutations $\xi$ of $\{0,\ldots,p-1\},$
\[
\bigl(W_{\xi(0)},\ldots,W_{\xi(p-1)}\bigr) \mbox{ has the same law as } (W_0,\ldots,W_{p-1}).
\]
Let $\overline{W}=p^{-1}\sum_{i=0}^{p-1}W_i$, for all $t$ in
$L^1(P)$, let
$P_A^W$ be the block-resampling empirical process defined by
\[
P_A^Wt=\frac1p\sum_{i=0}^{p-1}W_iL_qt(A_i).
\]
For all integrable random variables $F(X_1,\ldots,X_n,W_0,\ldots
,W_{p-1})$, let
\begin{eqnarray*}
&&\mathbb{E}_W[F(X_1,\ldots,X_n,W_0,\ldots,W_{p-1})]\\
&&\qquad=\mathbb{E}
[F(X_1,\ldots,X_n,W_0,\ldots,W_{p-1})|X_1,\ldots,X_n].
\end{eqnarray*}
Let $((\psi_{\lambda})_{\lambda\in\Lambda_m})_{m \in\mathcal
{M}_n}$ be
orthonormal bases of $(S_m)_{m \in\mathcal{M}_n}$ and let\break
$(\widehat{s}_{A,m}^W)_{m \in\mathcal{M}_n}$ be
the collection of resampling projection estimators
\[
\widehat{s}_{A,m}^W=\sum_{\lambda\in\Lambda_m}(P_A^W\psi_{\lambda
})\psi
_{\lambda}.
\]
The block-resampling penalties are defined as block-resampling
estimators of the ideal penalty by
%
%
\begin{equation}\label{defRespen}
\operatorname{pen}_W(m,C)=C\mathbb{E}_W\bigl(2(P_A^W-\overline
{W}P_A)(\widehat{s}^W_{A,m})\bigr).
\end{equation}
The idea of resampling is to mimic the behavior of the empirical
process~$P_A$ around~$P$ by the behavior of the resampling empirical
process $P_A^W$ around $\overline{W}P_{A}$. The resampling procedure
is a
plug-in method where the unknown functionals $F(P,P_n)$ are estimated
by $F(\overline{W}P_n, P_n^W)$. Hence, $\widehat{s}_{A,m}$ in
$\operatorname{pen}_{\mathrm{id}}(m)$ is replaced
by $\widehat{s}_{A,m}^W$ in $\operatorname{pen}_W(m,C)$ and, instead
of applying the process
$P_A-P$, we apply the process $P_A^W-\overline{W}P_A$. We take the expectation
with respect to the distribution of the resampling scheme to stabilize
the procedure. Finally, we let a normalizing constant $C$ free for this
general definition.\looseness=-1

We use a block-resampling scheme instead of a classical exchangeable
resampling scheme in order to preserve the dependence of the data
inside the blocks. This is a key point for the procedure to work.
Examples of resampling schemes can be found in \citet{Ar08}. The
classical block-bootstrap [\citet{Ku89}; \citet{LS92}] is obtained when the
distribution of $(W_0,\ldots,W_{p-1})$ is the multinomial $\mathcal{M}
(p,1/p,\ldots,1/p)$.
%
\subsection{The slope algorithm}
The ``slope heuristic'' has been introduced by \citet{BM07} in order to
calibrate the leading constant in a~penalty term [e.g., the
constant $C$ in (\ref{defRespen})]. It is based on the behavior of the
complexity of the selected model [recall the definition (\ref{defPPE})].
It states that there exist a family $(\Delta_m)_{m \in
\mathcal{M}_n}$ and
a constant $K_{\min}$ satisfying the following properties:
\begin{longlist}[(SH1)]
\item[(SH1)] When $\operatorname{pen}(m)\!\leq\!K\Delta_m$, with
$K\,{<}\,K_{\min}$, then $\Delta_{\widehat{m}}\!\geq\! c_1\max_{m \in\mathcal{M}_n}\Delta_m$.
\item[(SH2)] When $\operatorname{pen}(m)= K\Delta_m$, with $K>K_{\min
}$, then $\Delta
_{\widehat{m}}$ is much smaller.
\item[(SH3)] When $\operatorname{pen}(m)= 2K_{\min}\Delta_m$, then
$\widehat{s}_A$ satisfies a
sharp oracle inequality.
\end{longlist}
Based on this heuristic, \citet{BM07} introduced the following slope
algorithm. It can be used in practice when a family $(\Delta_m)_{m \in
\mathcal{M}_n}$
satisfying the slope heuristic is known.
\begin{itemize}
\item For all $K>0$, compute $\Delta_{\widehat{m}(K)}$ where $\widehat{m}(K)$
is defined as in (\ref{defPPE}) with $\operatorname{pen}(m)=K\Delta_m$.
\item Find $\tilde{K}$ such that $\Delta_{\widehat{m}(K)}$ is very large
for $K<\tilde{K}$ and much smaller when $K>\tilde{K}$.
\item Choose the final $\widehat{m}$ equal to $\widehat{m}(2\tilde{K})$.
\end{itemize}
The idea is that $\tilde{K}\sim K_{\min}$ since we observe a jump of
the complexity of the selected model around $K=\tilde{K}$ [thanks to
(SH1), (SH2)] and thus that the final estimator, selected by the
penalty $2\tilde{K}\Delta_m\sim2K_{\min}\Delta_m$, satisfies an
optimal oracle inequality [by (SH3)].
%
%
\subsection{Some measures of dependence}\label{Ch4ss23}
\subsubsection{\texorpdfstring{$\beta$-mixing data}{beta-mixing data}}\label{Ch4ss231}
Volkonski{\u\i} and Rozanov
(\citeyear{RV59}) defined the coefficient~$\beta$ as follows. Let $Y$ be a
random variable defined on a probability space $(\Omega,\mathcal
{A},\mathbb{P}
)$ and let $\mathcal{M}$ be a $\sigma$-algebra in $\mathcal{A}$, let
\[
\beta(\mathcal{M},\sigma(Y))=\mathbb{E}\Bigl(\sup_{A\in\mathcal
{B}}|\mathbb{P}_{Y|\mathcal{M}}(A)-\mathbb{P}
_Y(A)|\Bigr).
\]
For all stationary sequences of random variables $(X_n)_{n\in\mathbb{Z}}$
defined on $(\Omega,\mathcal{A},\mathbb{P})$, let
\[
\beta_k=\beta\bigl(\sigma(X_i, i\leq0),\sigma(X_i, i\geq k)\bigr).
\]
The process $(X_n)_{n\in\mathbb{Z}}$ is said to be $\beta$-mixing
when $\beta
_k\rightarrow0$ as $k\rightarrow\infty$.
%
\subsubsection{\texorpdfstring{$\tau$-mixing data}{tau-mixing data}}\label{Ch4ss232}
Dedecker and Prieur
(\citeyear{DP05}) defined the coefficient~$\tau$ as follows. For all $l$ in
$\mathbb{N}^*$, for all $x,y$ in $\mathbb{R}^l$, let $d_l(x,y)=\sum
_{i=1}^l|x_i-y_i|$.
For all $l$ in~$\mathbb{N}^*$, for all functions~$t$ defined on
$\mathbb{R}^l$, the Lipschitz semi-norm of~$t$ is defined by
\[
\operatorname{Lip}_{l}(t)=\sup_{x\neq y\in\mathbb{R}^l}\frac
{|t(x)-t(y)|}{d_l(x,y)}.
\]
For all functions $t$ defined on $\mathbb{R}$, we will denote for
short by $\operatorname{Lip}
(t)=\operatorname{Lip}_1(t)$. Let~$\lambda_1$ be the set of all
functions $t\dvtx \mathbb{R}
^l\rightarrow\mathbb{R}$ such that $\operatorname{Lip}_{l}(t)\leq
1$. For all integrable, $\mathbb{R}
^l$-valued, random variables $Y$ defined on a probability space
$(\Omega
,\mathcal{A},\mathbb{P})$ and all $\sigma$-algebra $\mathcal{M}$ in
$\mathcal{A}$, let
\[
\tau(\mathcal{M},Y)=\mathbb{E}\Bigl(\sup_{t\in\lambda_1}|\mathbb
{P}_{Y|\mathcal{M}}(t)-\mathbb{P}_Y(t)|\Bigr).
\]
For all stationary sequences of integrable random variables
$(X_n)_{n\in
\mathbb{Z}}$ defined on $(\Omega,\mathcal{A},\mathbb{P})$, for all
integers $k,r$, let
\[
\tau_{k,r}=\max_{1\leq l\leq r}\frac{1}{l}\sup_{k\leq i_1<\cdots<
i_l}\bigl\lbrace\tau\bigl(\sigma(X_p, p\leq0),(X_{i_1},\ldots
,X_{i_l})\bigr)\bigr\rbrace,
\qquad\tau_k=\sup_{r\in\mathbb{N}^*}\tau_{k,r}.
\]
The process $(X_n)_{n\in\mathbb{Z}}$ is said to be $\tau$-mixing
when $\tau
_k\rightarrow0$ as $k\rightarrow\infty$.

\subsection{Main assumptions}
%
\subsubsection{\texorpdfstring{A specific collection for $\tau$-mixing sequences}{A specific collection for tau-mixing sequences}}
Wavelet spaces have been widely used in density estimation since the
oracle is adaptive over Besov spaces [see \citet{BM97}].

\textit{Dyadic wavelet spaces:} Let $r$ be a real number, $r\geq1$. We
work with an $r$-regular orthonormal multiresolution analysis of
$L^2(\mu)$, associated with a~compactly supported scaling function
$\phi
$ and a compactly supported mother wavelet~$\psi$. Without loss of
generality, we suppose that the support of the functions $\phi$ and~$\psi$ is
included in an interval $[A_1,A_2)$ where $A_1$ and $A_2$ are
integers such that $A_2-A_1=A\geq1$. For all $k$ in $\mathbb{Z}$ and
$j$ in $\mathbb{N}
^*$, let $\psi_{0,k}\dvtx  x\rightarrow\sqrt{2}\phi(2x-k)$ and $\psi
_{j,k}\dvtx  x\rightarrow2^{j/2}\psi(2^jx-k)$. The family $\{(\psi
_{j,k})_{j\geq0,k\in\mathbb{Z}}\}$ is an orthonormal basis of
$L^2(\mu
)$. The collection of dyadic wavelet spaces is described as follows.

\begin{longlist}[W]
\item[{[W]}] \textit{dyadic wavelet generated spaces: let $J_n=[\log
_2(n)]$, for all $J_m=1,\ldots,\allowbreak J_n$, let}
\[
\Lambda_m=\{(j,k), 0\leq j\leq J_m, k\in\mathbb{Z}\}
\]
\textit{and let $S_m$ be the linear span of $\{\psi_{\lambda}\}_{\lambda\in
\Lambda_m}$.}\vspace*{8pt}
\end{longlist}

\subsubsection{General framework} We present in this section a set of
assumptions sufficient to prove the theorems. None of them is used to
build the penalties.

\begin{longlist}[(H1)]
\item[(H1)] \textit{There exists a constant $\kappa_a$ such that,
for all $m$, $m'$ in $\mathcal{M}_n$, for all~$t$ in $S_m+S_{m'}$, with $\|t\|\leq1$,
there exist $t_m$ in $S_m$ and $t_{m'}$ in $S_{m'}$, with $\|t_m\|\vee\|t_{m'}\|\leq\kappa_a$
such that $t=t_m+t_{m'}$.}
\end{longlist}

(H1) is typically satisfied for nested collections as [W].
\begin{longlist}[(H2)]
\item[(H2)] \textit{$N_n=\operatorname{Card}(\mathcal{M}_n)$ is finite and
there exist
constants $c_\mathcal{M}$, $\alpha_\mathcal{M}$ such that $N_n\leq
c_\mathcal{M}n^{\alpha_\mathcal{M}}$.}
\end{longlist}

(H2) means that the collection is not too rich and thus that the
model selection problem is not too hard. It is satisfied by the
collection [W].

Let us introduce some notation. For all $m$ in $\mathcal{M}_n$, for
all orthonormal bases $(\psi_{\lambda})_{\lambda\in\Lambda_m}$ of
$S_m$, let
\begin{eqnarray*}
D_{A,m}&=&q\sum_{\lambda\in\Lambda_m}\operatorname{Var}(L_q(\psi
_{\lambda})(A_0)),\qquad
R_{A,m}=n\|s-s_m\|^2+2D_{A,m},\\
B_m&=&\{t\in S_m,\|t\|\leq1\}, \qquad b_{m}=\sup_{t\in B_m}\|t\|
_{\infty}.
\end{eqnarray*}
$D_{A,m}$, and thus $R_{A,m}$, are well defined since we can check with
Cauchy--Schwarz inequality that
\[
D_{A,m}=q\mathbb{E}\Bigl[\Bigl(\sup_{t\in B_m}L_qt(A_0)-Pt\Bigr)^2\Bigr].
\]
Two quantities will play a fundamental role to discuss the results. The
first one is the risk of an oracle:
\[
R_n=\inf_{m \in\mathcal{M}_n}R_{A,m}.
\]
We are typically interested in non parametric problems where $R_n/n\sim
n^{-\gamma}$ for some $0<\gamma<1$. This situation occurs, for example,
when $s$ is a regular function, in this case, we have $R_n/n= \kappa
n^{-2\alpha/(2\alpha+1)}$, for some $\alpha>0$, $\kappa>0$. We will
make the following assumption:

\begin{longlist}[(H3)]
\item[(H3)] \textit{There exists a constant $\kappa_R>0$ such that
$R_n\geq
\kappa_R(\ln n)^8$.}
\end{longlist}

The constant $8$ in $(\ln n)^{8}$ is technical, it yields the rate
$\varepsilon_n=(\ln n)^{-1/2}$ in the oracle inequalities. \citet{Ar08}
replaced this assumption by a lower bound on the bias of the models. It
implies that $R_n\geq\kappa n^{\gamma}$, for some constants $\kappa
>0$, $1>\gamma>0$ and therefore assumption (H3).

\begin{longlist}[(H4)]
\item[(H4)] \textit{There exists a constant $c_D>0$ such that}
\[
\forall m\in\mathcal{M}_n\qquad P\Bigl( \sup_{t\in B_{m}}t^{2} \Bigr)\geq c_Db_m^2.
\]
\end{longlist}

It is shown in the \hyperref[SectionAppendix]{Appendix} that some classical
examples of collections $(S_{m}, m \in\mathcal{M}_n)$ as regular
histograms, Fourier
spaces and [W] satisfy (H4).

%
The following assumptions will be used to prove the slope heuristic. We
introduce a second quantity, that will play a fundamental role. Let
\[
D_n^*=\max_{m \in\mathcal{M}_n}D_{A,m}.
\]
In classical collection of models, like [W], $D_n^*\sim c n$.
This is why we introduce the following assumption:

\begin{longlist}[(H5)]
\item[(H5)] \textit{$D_n^*/R_n\rightarrow\infty$ when $n$ grows
to infinity.}
\end{longlist}

We will prove that, when the data are mixing. $D_{A,m}\simeq
n\mathbb{E}(\|s_m-\widehat{s}_{A,m}\|^2)$ represents the variance
term of
the risk. It is a natural measure of the complexity of the models.
Hence, $D_n^*$ represents the maximal complexity of the models.
Moreover, $R_n$ is the risk of the oracle. It balances the complexity
and the bias term and has therefore the same order as the complexity of
an oracle. Hence, assumption (H5) means that the largest complexity
in the collection $(S_m)_{m \in\mathcal{M}_n}$ is much larger than
the one of an oracle, which is a~na\-tural condition for the slope heuristic to hold.
We need a final assumption.\looseness=-1

\begin{longlist}[(H6)]
\item[(H6)] \textit{For all $m^*$ such that $D_{A,m^*}=D^*_n$, we have}
\[
\frac{n\|s-s_{m^*}\|^2}{D_n^*}\rightarrow0 \qquad\mbox{when } n\rightarrow
\infty.
\]
\end{longlist}

When $D_n^*$ is of order $n$, (H6) simply means that the distance
between $s$ and a complex model goes to $0$. In general, it means that
for these complex models, the bias part of the risk is negligible
compared to the variance part. We conclude this section by the
assumptions on the mixing coefficients. All mean that these
coefficients are sufficiently small. Let $\gamma=\beta$ or $\tau$.

\begin{longlist}[$\mathrm{S}(\beta)$]
\item[{[$\operatorname{AR}(\theta)$]}] \textit{arithmetical $\gamma$-mixing with
rate $\theta$: there exists $C>0$ such that, for all $k$ in $\mathbb{N}$,
$\gamma_k\leq C(1+k)^{-(1+\theta)}$.}

\item[$\mathrm{S}(\beta)$] \textit{$\sum_{l\geq1}(l+1)\beta_{l}\leq c_D/64$,
where $c_D$ is defined in} (H4).

We prove in the \hyperref[SectionAppendix]{Appendix} that $c_{D}=1$ for regular histograms and Fourier spaces.

\item[$\mathrm{S}(\tau, \mathrm{W})$] \textit{$\sum_{l\geq1}( {s}^2\tau _{l}
)^{1/3}\leq C(W)$, where $C(W)$ depends only on $\phi$, $\psi$.}

The value of the constant $C(W)$ is given in Lemma \ref{lemconcBeta} of~\citet{le11}.
\end{longlist}

\section{\texorpdfstring{Results for $\tau$-mixing sequences}{Results for tau-mixing sequences}}\label{Ch4S3}
\subsection{Resampling penalties}\vspace*{3pt}
The result of this section is that PPE selected by block-resampling
penalties satisfy sharp oracle inequalities.
\begin{theo}\label{theoRestau}
Let $X_1,\ldots,X_n$ be a strictly stationary sequence of real valued
random variables with common density $s$ and let\vadjust{\goodbreak} $(S_m)_{m \in\mathcal
{M}_n}$ be a~collection of regular wavelet spaces [\textup{W}] satisfying (\textup{H3}), (\textup{H4}).
Let $p$, $q$ be two integers such that $2pq=n$ and $\frac12\sqrt
{n}(\ln n)^2\leq p\leq\sqrt{n}(\ln n)^2.$

Let $\tilde{C}_{W}=\operatorname{Var}(W_1-\overline{W})^{-1}$,
$C>\tilde{C}_{W}/2$ and let $\tilde{s}_A$ be the
PPE defined in~(\ref{defPPE}) with the penalty $\operatorname
{pen}_W(m,C)$ defined in~(\ref{defRespen}).

Assume that there exists $\theta>5$ such that $X_1,\ldots,X_n$ are
arithmetically [$\operatorname{AR}(\theta)$] $\tau$-mixing and satisfy
$\mathrm{S}(\tau,\mathrm{W})$. Let $\varepsilon_{n}=(\ln n)^{-1/2}$,
$\kappa(C)=\break
2(C\tilde{C}_{W}^{-1}-1 )$.

There exist constants $\kappa_1$, $\kappa_2$ such that we have
%
%
\begin{equation}\label{OEMR}
K_n\mathbb{E}(\|s-\tilde{s}_A\|^2)\leq\mathbb{E}\Bigl(\inf_{m \in
\mathcal{M}_n}\|s-\widehat{s}
_{A,m}\|^2\Bigr)+\frac{\kappa_2}{n},
\end{equation}
with
\[
K_n=\frac{(1\wedge(1+\kappa(C))-\kappa_{1}\varepsilon_{n}}{(1\vee
(1+\kappa(C))+\kappa_{1}\varepsilon_{n}}.
\]
\end{theo}
\textit{Comments:}
\begin{itemize}
%
\item The constant $C$ has to be chosen asymptotically equal to $\tilde
{C}_{W}$.
If we choose $C>\tilde{C}_{W}$, we still get an oracle inequality,
with a leading
constant less sharp. On the other hand, if we choose $C<\tilde{C}_{W}$
we can
have $K_{n}\leq0$ in~(\ref{OEMR}). This is a first reason why it is
generally useful to over-penalize a~little bit from a~nonasymptotic
point of view.
\end{itemize}
%
\vspace*{3pt}\subsection{Slope heuristic}
Theorem \ref{theoRestau} gives a totally data driven penalty which
satisfies a sharp oracle inequality, therefore, the heuristic is not
necessary to obtain asymptotically optimal results. However, $C$ can be
optimized for small samples. Moreover, the slope algorithm is faster to
compute than resampling penalties when a deterministic quantity can be
used in the slope heuristic. Theorem \ref{theoSlope1tau} hereafter
justifies property (SH1) of the heuristic. $\Delta_m$ is the
variance term $D_{A,m}/n$ and $K_{\min}=2$.
\begin{theo}\label{theoSlope1tau}
Let $X_1,\ldots,X_n$ be a strictly stationary sequence of real valued
random variables with common density $s$ and let $(S_m)_{m \in\mathcal
{M}_n}$ be a~collection of regular wavelet spaces [\textup{W}] satisfying (\textup{H3})--(\textup{H6}).
Let $p$, $q$ be two integers such that $2pq=n$ and $\frac
12\sqrt{n}(\ln n)^2\leq p\leq\sqrt{n}(\ln n)^2$.

Assume that there exists a constant $0<\delta<1$ such that, for all
$m$ in $\mathcal{M}_n$,
%
%
\begin{equation}\label{condpentropetite}
0\leq\operatorname{pen}(m)\leq(2-\delta)\frac{D_{A,m}}n,
\end{equation}
and let $\tilde{s}_A$ be the $\mathit{PPE}$ defined in (\ref{defPPE}).

Assume that there exists $\theta>5$ such that $X_1,\ldots,X_n$ are
arithmetically [$\operatorname{AR}(\theta)$] $\tau$-mixing and satisfy $\mathrm{S}(\tau,\mathrm{W})$.
There exist constants $\kappa_1$, $\kappa_2$ such that
%
%
\begin{equation}\label{SME}
\mathbb{E}(D_{A,\widehat{m}})\geq\frac{4\delta}{9}D^*_{n}-\kappa_1.
\end{equation}
%
%
\begin{equation}\label{SME2}
\mathbb{E}(\|s-\tilde{s}_{A}\|^2)\geq\frac{\delta}{5}\frac
{D_n^*}{R_n}\biggl(\mathbb{E}\Bigl(\inf_{m \in\mathcal{M}_n}\|s-\widehat
{s}_{A,m}\|^2\Bigr)-\frac
{\kappa_2}n\biggr).\vspace*{3pt}
\end{equation}
\end{theo}
\textit{Comments:}
\begin{itemize}
\item Inequality (\ref{SME}) states that $D_{A,\widehat{m}}$ is as large as
possible when the penalty term is too small. This is exactly (SH1)
with $\Delta_m=D_{A,m}$.
\item Inequality (\ref{SME2}) states that the model selected by a too
small penalty is never an oracle. This is another reason why it is
interesting to choose $C>\tilde{C}_{W}$ in Theorem \ref{theoRestau}.
%
\end{itemize}
The following theorem justifies properties (SH2), (SH3) of the slope
heuristic.
\begin{theo}\label{theoSlope2tau}
Let $X_1,\ldots,X_n$ be a strictly stationary sequence of real valued
random variables with common density $s$ and let $(S_m)_{m \in\mathcal
{M}_n}$ be a~collection of regular wavelet spaces [\textup{W}] satisfying (\textup{H3}), (\textup{H4}).
Let $p$, $q$ be two integers such that $2pq=n$ and $\frac12\sqrt
{n}(\ln n)^2\leq p\leq\sqrt{n}(\ln n)^2.$

Assume that there exist $\delta_+\geq-\delta_->-1$, $\varepsilon
\geq0$
and some constants $\kappa_1$, $\kappa_2$ satisfying, for all $m$ in
$\mathcal{M}_n$,
%
\begin{eqnarray}
\mathbb{E}\biggl[\sup_{m \in\mathcal{M}_n}\biggl((2-\delta_{-})\frac
{2D_{A,m}}{n}-\operatorname{pen}
(m)-\varepsilon\frac{R_{A,m}}n\biggr)_+\biggr]&\leq&\frac{\kappa
_1}n,\label{condpenassezgrande}\\
\mathbb{E}\biggl[\sup_{m \in\mathcal{M}_n}\biggl(\operatorname
{pen}(m)-(2+\delta_{+})\frac
{D_{A,m}}{n}-\varepsilon\frac{R_{A,m}}n\biggr)_+\biggr]&\leq&\frac
{\kappa_2}n.\label{condpenpastropgrande}
\end{eqnarray}
Let $\tilde{s}_A$ be the \textit{PPE} defined in (\ref{defPPE}) with
$\operatorname{pen}(m)$
and let $\varepsilon_n=(\ln n)^{-1/2}$.

Assume that there exists $\theta>5$ such that $X_1,\ldots,X_n$ are
arithmetically [$\operatorname{AR}(\theta)$] $\tau$-mixing and satisfy $\mathrm{S}(\tau,\mathrm{W})$.
There exist constants $\kappa_1$, $\kappa_2$, $\kappa_3$,
such that
%
%
\begin{equation}\label{eqOISlopetau}
K_n\mathbb{E}(\|\tilde{s}_A-s\|^2)\leq\mathbb{E}\Bigl(\inf_{m \in
\mathcal{M}_n}\|s-\widehat{s}
_{A,m}\|^2\Bigr)+\frac{\kappa_2}n,
\end{equation}
with
\[
K_n=\frac{(1\wedge( 1-\delta_- ))-\kappa_1(\varepsilon
_n+\varepsilon
)}{(1\vee(1+\delta_+))+\kappa_1(\varepsilon_n+ \varepsilon)}.
\]
Moreover, we have
%
%
\begin{equation}\label{eqDhmpetit}
K_n\mathbb{E}(D_{A,\widehat{m}})\leq R_n+\kappa_3.
\end{equation}

\end{theo}
\textit{Comments:}
\begin{itemize}
\item When $\operatorname{pen}(m)$ becomes larger than $2D_{A,m}/n$,
$D_{A,\widehat{m}}$
jumps from $D_n^*$ (\ref{SME}) to $R_n$ [(\ref{eqDhmpetit}) for
$\delta
_+$ and $-\delta_-$ close to $-1$]. This justifies (SH2) since
$R_n\ll D_n^*$.
\item A model selected with a penalty $4D_{A,m}/n$ satisfies an oracle
inequality (Theorem \ref{theoSlope2tau} for $\delta_+$ and $\delta_-$
close to $0$). This justifies (SH3).
\item $D_{A,m}$ is unknown and cannot be used in the slope algorithm.
We show [Lemma~\ref{lemconcBeta} in \citet{le11}] that $D_{A,m}$ satisfies $\kappa
_{*}2^{J_{m}}\leq D_{A,m}\leq\kappa^{*}2^{J_{m}}$. The slope heuristic
might hold for $\Delta_m=2^{J_m}/n$, but a complete proof requires
moreover that $\kappa_{*}\simeq\kappa^{*}$. However, we obtain in the
proof of Theorem \ref{theoRestau} that $\operatorname{pen}_W(m,
\tilde{C}_{W})$ satisfies
(\ref{condpenassezgrande}) and (\ref{condpenpastropgrande}) for
$\delta_+=\delta_-=0$ and $\varepsilon=\kappa\varepsilon_n$. Since
(\ref{condpentropetite}) can be modified to work with random penalties, we
can apply the slope algorithm with $\operatorname{pen}_W(m, 1)$
instead of $D_{A,m}/n$.

\end{itemize}

\section{\texorpdfstring{Results for $\beta$-mixing sequences}{Results for beta-mixing sequences}}\label{Ch4S4}
We show that block-resampling penalties select oracles and that the
slope heuristic holds in this case.
%
\subsection{Resampling penalties}
\begin{theo}\label{theoResbeta}
Let $X_1,\ldots,X_n$ be a strictly stationary sequence of real valued
random variables with common density $s$ and let $(S_m)_{m \in\mathcal
{M}_n}$ be a~collection of linear spaces satisfying (\textup{H1})--(\textup{H4}). Let $p$,
$q$ be two integers such that $2pq=n$ and $\frac12\sqrt{n}(\ln
n)^2\leq p\leq\sqrt{n}(\ln n)^2.$

Let $\tilde{C}_{W}=\operatorname{Var}(W_1-\overline{W})^{-1}$,
$C>\tilde{C}_{W}/2$ and let $\tilde{s}_A$ be the
\textit{PPE} defined in~(\ref{defPPE}) with the block-resampling penalty
$\operatorname{pen}_W(m,C)$ defined in~(\ref{defRespen}).

Assume that there exists $\theta>2$ such that $X_1,\ldots,X_n$ are
arithmetically [$\operatorname{AR}(\theta)$] $\beta$-mixing and satisfy $\mathrm{S}(\beta)$.
Let $\varepsilon_{n}=(\ln n)^{-1/2}, \kappa(C)=2(
C\tilde{C}_{W}^{-1}-1 ).$

There exist constants $\kappa_1$, $\kappa_2$ such that
%
%
\begin{equation}\label{OPMR}
\qquad P\Bigl( K_n\|s-\tilde{s}_A\|^2\leq\inf_{m \in\mathcal{M}_n}\|s-\widehat{s}_{A,m}\|^2 \Bigr)\geq 1-
\kappa_2\biggl(\frac{1}{n^2}\vee\frac{(\ln n)^{4+2\theta}}{n^{\theta/2}}\biggr),
\end{equation}
with
\[
K_n=\frac{(1\wedge(1+\kappa(C)))-\kappa_{1}\varepsilon_{n}}{(1\vee
(1+\kappa(C))+\kappa_{1}\varepsilon_{n}}.
\]
\end{theo}
\textit{Comments:}
\begin{itemize}
\item The coupling lemma of \citet{Be79} for $\beta$-mixing processes
is much stronger than the one satisfied by $\tau$-mixing data [\citet{DP05}]. This is why Theorem \ref{theoResbeta} covers more
collections of models than Theorem \ref{theoRestau} and why we prove
oracle inequalities in probability.
\end{itemize}\vspace*{-5pt}
%
\subsection{Slope heuristic} The following theorems are adaptations to
the $\beta$-mixing case of Theorems \ref{theoSlope1tau} and \ref{theoSlope2tau}.\vspace*{-3pt}
\begin{theo}\label{theoSlope1beta}
Let $X_1,\ldots,X_n$ be a strictly stationary sequence of real valued
random variables with common density $s$ and let $(S_m)_{m \in\mathcal
{M}_n}$ be a~collection of linear spaces satisfying (\textup{H1})--(\textup{H6}).
Let $p$, $q$ be two integers such that $2pq=n$ and $\frac12\sqrt
{n}(\ln n)^2\leq p\leq\sqrt{n}(\ln n)^2.$

Let $\tilde{s}_A$ be the \textit{PPE} defined in (\ref{defPPE}) with a penalty
$\operatorname{pen}(m)$ satisfying, for all $m$ in $\mathcal{M}_n$,
condition (\ref{condpentropetite}) of Theorem \ref{theoSlope1tau}.

Assume that there exists $\theta>2$ such that $X_1,\ldots,X_n$ are
arithmetically [$\operatorname{AR}(\theta)$] $\beta$-mixing and satisfy $\mathrm{S}(\beta)$.
There exists a constant $\kappa$ and an event~$\Omega_n$ such that
\[
\mathbb{P}(\Omega_n)\geq1-\kappa\biggl(\frac{1}{n^2}\vee\frac{(\ln
n)^{4+2\theta
}}{n^{\theta/2}}\biggr),
\]
and, on $\Omega_n$,
%
%
\begin{equation}\label{eqpasdOI}
D_{A,\widehat{m}}\geq\frac{4\delta}{9}D^*_{n},\qquad \|s-\tilde{s}_{A}\|
^2\geq
\frac{\delta}{5}\frac{D_n^*}{R_n}\inf_{m \in\mathcal{M}_n}\|
s-\widehat{s}_{A,m}\|^2.\vspace*{-3pt}
\end{equation}
\end{theo}
\begin{theo}\label{theoSlope2beta}
Let $X_1,\ldots,X_n$ be a strictly stationary sequence of real valued
random variables with common density $s$ and let $(S_m)_{m \in\mathcal
{M}_n}$ be a~collection of linear spaces satisfying (\textup{H1})--(\textup{H4}). Let $p$,
$q$ be two integers such that $2pq=n$ and $\frac12\sqrt{n}(\ln
n)^2\leq p\leq\sqrt{n}(\ln n)^2.$

Assume that there exist $\delta_+\,{\geq}\,-\delta_-\,{>}\,-1$,
$\varepsilon\,{\geq}\,0$, $0\,{\leq}\,\eta\,{<}\,1$ and an event~$\Omega_{\operatorname{pen}}$, with
$\mathbb{P}(\Omega_{\operatorname{pen}
})\geq1-\eta$ such that, on~$\Omega_{\operatorname{pen}}$, for all
$m$ in $\mathcal{M}_n$,
%
%
\begin{equation}\label{condbonnePen}
(2-\delta_{-})\frac{2D_{A,m}}{n}-\varepsilon\frac{R_{A,m}}n\leq
\operatorname{pen}
(m)\leq(2+\delta_{+})\frac{2D_{A,m}}{n}+\varepsilon\frac{R_{A,m}}n.
\end{equation}
Let $\tilde{s}_A$ be the \textit{PPE} defined in (\ref{defPPE}) with
$\operatorname{pen}$.

Assume that there exists $\theta>2$ such that $X_1,\ldots,X_n$ are
arithmetically [$\operatorname{AR}(\theta)$] $\beta$-mixing and satisfy $\mathrm{S}(\beta)$.
There exist constants $\kappa_1$, $\kappa_2$ and an event
$\Omega_n^*$ such that
\[
\mathbb{P}(\Omega_n^*)\geq1-\eta-\kappa_2\biggl(\frac{1}{n^2}\vee\frac
{(\ln
n)^{4+2\theta}}{n^{\theta/2}}\biggr),
\]
and, on $\Omega_n^*$,
%
%
\begin{equation}\label{eqOISlopebeta}
K_n\|\tilde{s}_A-s\|^2\leq\inf_{m \in\mathcal{M}_n}\|s-\widehat
{s}_{A,m}\|^2,
\end{equation}
with
\[
K_n= \frac{(1\wedge( 1-\delta_- ))-\kappa_1(\varepsilon_n+
\varepsilon)}{(1\vee(1+\delta_+))+\kappa_1(\varepsilon_n+
\varepsilon)}.
\]
Moreover, $\Omega_n^*$, $2K_nD_{A,\widehat{m}}\leq3R_n.$
\end{theo}
\textit{Comments:}
\begin{itemize}
\item We refer to the comments of Theorems \ref{theoSlope1tau} and
\ref{theoSlope2tau} where we explain why Theorems \ref{theoSlope1beta} and \ref{theoSlope2beta} imply the slope heuristic
with $\Delta_m=D_{A,m}/n$, $K_{\min}=2$.
\item As in Theorem \ref{theoSlope2tau}, $D_{A,m}$ cannot be used to
build a model selection procedure. A deterministic shape of $D_{A,m}$
is unknown, although we prove in the Supplementary Material that
$D_{A,m}$ is bounded by $b_m^*$. However, $\operatorname{pen}_W(m,1)$
can be used
instead of $D_{A,m}$.
\end{itemize}
%
%
\subsection{Discussion and perspectives}\label{Ch4ss43}
Block-resampling penalties yield data driven procedures for the
estimation of the marginal density in a mixing framework. The selected
estimators satisfy sharp oracle inequalities without remainder term.
This improves Theorems \ref{theoRestau} and \ref{theoResbeta} in \citet{Le08} and Theorem \ref{theoRestau} in
\citet{CM02}, where the leading constants was built with the mixing
coefficients of the process. Moreover, our results hold for possibly
infinite dimensional models.

\citet{La08} gave also a model selection procedure to estimate the
stationary density and the transition probability of a Markov Chain.
She worked with a stationary chain, irreducible, aperiodic and
positively recurrent, which is therefore $\beta$-mixing. Her density
estimator is selected by a~penalty equal to $Kd_m/n$ with a constant
$K$ that ``depends on the law of the chain'' [see Remark 4 after Theorem~3 in \citet{La08}]. She proposed to estimate $K$ in the simulations by
the slope algorithm. We prove the slope heuristic, justifying that the
slope algorithm can be used to optimize the leading constant. It would
be interesting to see if resampling penalties may be used in her
context to estimate the transition probabilities.

\citet{GW09} worked with other weak mixing coefficients [namely
$\lambda
$ and $\tilde{\phi}$; see \citet{DDLLLP} for a definition] and
studied a
wavelet thresholded estimator. The main advantage is that the
thresholded estimator is adaptive over a larger class of Besov spaces
than the oracle over the collection [W] [for details about this
important issue see \citet{BBM99}]. The main drawback is that their
threshold is built with the mixing coefficients.

Block-resampling penalties can be extended to the statistical learning
framework of \citet{MN06}, where the slope algorithm has already been
defined [\citet{AM08}]. We believe that these procedures perform well in
this context but the problem remains open.

The main drawback of our approach is that we use only $n/2$ data.
Moreover, the deterministic choice of the number $p$ of blocks is not
optimized. For example, when the data are geometrically $\beta$-mixing,
which means that, for some constants $\theta>0$, $C>0$, $\beta_k\leq
Ce^{-\theta k}$, choosing $p$ of order $n(\ln n)^{-2}$ would improve
the rates of convergence of the leading constant. An interesting
direction of research would be to provide data-driven choices of $p$
and $q$ to improve these rates, and a data-driven choice of blocks to
use more data.

In practice, the computation time is also a very important issue.
Actually, the conditional expectation is a bit long to evaluate and
some efforts have to be done in this direction. Things can be improved
if we obtain a~deterministic shape of the ideal penalty, as in the
independent case, since the slope heuristic is faster to compute with a
deterministic $\Delta_m$. We obtain upper and lower bounds on
$\operatorname{pen}_{\mathrm{id}}$,
but our inequalities are not sharp enough to justify completely the
slope heuristic. We can also think of the $V$-fold cross validation
penalties defined in \citet{ArVF}. These penalties are also faster to
compute than the resampling penalties. They can be viewed as resampling
penalties defined with nonexchangeable weights. These issues are far
beyond the objectives of the present paper and will be addressed in
forthcoming works.

\vspace*{-3pt}\section{Proofs}\vspace*{-3pt}
\label{Ch4S5}
\subsection{Notation}\label{Ch4ss51}
Recall that $p$ and $q$ are integers such that $2pq=n$, and that $\sqrt
{n}(\ln n)^2/2\leq p\leq\sqrt{n}(\ln n)^2$. For all $k=0,\ldots
,p-1$, let
$I_k=(2kq+1,\ldots,(2k+1)q)$, $A_k=(X_i)_{i\in I_k}$ and $I=\bigcup
_{k=0}^{p-1}I_k$. For all $t$ in $L^2(\mu)$ and all $x_1,\ldots,x_q$ in
$\mathbb{R}$,
\begin{eqnarray*}
L_q(t)(x_1,\ldots,x_q)&=&\frac1q\sum_{i=1}^qt(x_i),\qquad P_At=\frac1p\sum
_{k=0}^{p-1}L_q(t)(A_k)=\frac2n\sum_{i\in I}t(X_i),
\\
\nu_A(t)&=&(P_A-P)(t).
\end{eqnarray*}
For all $m$ in $\mathcal{M}_n$, we denote by $(\psi_{\lambda
})_{\lambda\in
\Lambda_m}$ an orthonormal basis of $S_m$. The estimator $\widehat{s}_{A,m}$
associated to the model $S_m$, is defined as
\[
\widehat{s}_{A,m}=\sum_{\lambda\in\Lambda_m}(P_A\psi_{\lambda
})\psi_{\lambda}.
\]
Classical computations show that, if $s_m$ denotes the orthogonal
projection of $s$ onto~$S_m$,
\[
s_m=\sum_{\lambda\in\Lambda_m}(P\psi_\lambda)\psi_\lambda,\qquad \mbox{hence }
\|\widehat{s}_{A,m}-s_m\|^2=\sum_{\lambda\in\Lambda_m}(\nu_A\psi_\lambda)^2.
\]
The ideal penalty, $2\nu_A(\widehat{s}_{A,m})$ satisfies
\[
\nu_A(\widehat{s}_{A,m}-s_m)+\nu_A(s_m)=\sum_{\lambda\in\Lambda
_m}(\nu_A\psi_\lambda)^2+\nu_A(s_m)=\|\widehat{s}_{A,m}-s_m\|^2+\nu_A(s_m).
\]
For all $m$, $m'$ in $\mathcal{M}_n$, let
\begin{eqnarray*}
p(m)&=&\|s_m-\widehat{s}_{A,m}\|^2=\sup_{t\in B_m}(\nu_A(t))^2=\sum
_{\lambda\in
\Lambda_m}(\nu_A(\psi_{\lambda}))^2,
\\
\delta(m,m')&=& 2\nu_A(s_m-s_{m'}).
\end{eqnarray*}
Hereafter $W_0,\ldots,W_{p-1}$ denotes a resampling scheme, $\overline{W}
=p^{-1}\sum_{i=0}^{p-1}W_i$, $P_A^{W}$ denotes the resampling empirical
process, defined for all measurable functions $t$ by
\[
P_A^Wt=\frac1p\sum_{i=0}^{p-1}W_iL_qt(A_i).
\]
We introduce also $\nu_A^W=P_A^W-\overline{W}P_A$ and $\tilde
{C}_W=(\operatorname{Var}
(W_1-\overline{W}))^{-1}$. For any orthonormal basis $(\psi_{\lambda
})_{\lambda
\in\Lambda_m}$ of $S_m$, let
\[
p_W(m)=\tilde{C}_W\sum_{\lambda\in\Lambda_m}\mathbb{E}_W((\nu
^W_A(\psi
_{\lambda}))^2).
\]
$p_W(m)$ is well defined since, from the Cauchy--Schwarz inequality,
\[
p_W(m)=\tilde{C}_W\mathbb{E}_W\Bigl(\sup_{t\in B_m}(\nu_A^Wt)^2\Bigr).
\]
Let $\varepsilon_{n}=(\ln n)^{-1/2}$ and let $\kappa>0$. Let
$\mathcal{M}$
denote one of the set $\mathcal{M}_{n}$ or $\mathcal{M}_{n}^{2}$.
When $\mathcal{M}=\mathcal{M}_{n}$, for
all $\overline{m}$ in $\mathcal{M}$ let $R_{A,\overline{m}}=R_{A,m}$
and when $\mathcal{M}=\mathcal{M}_{n}^{2}$,
for all $\overline{m}=(m,m')$ in $\mathcal{M}$, let $R_{A,\overline
{m}}=R_{A,m}\vee R_{A,m'}$.
For all $\overline{m}$ in $\mathcal{M}$, let
%
\begin{eqnarray}
f_1(\overline{m}, \kappa)&=&p(m)-\frac{2D_{A,m}}n-\kappa\varepsilon
_{n}\frac
{R_{A,m}}n,\label{deff1}\\
f_2(\overline{m},\kappa)&=&\frac{2D_{A,m}}n-p(m)-\kappa\varepsilon
_{n}\frac
{R_{A,m}}n,\label{deff2}\\
f_3(\overline{m},\kappa)&=&p(m)-p_W(m)-\kappa\varepsilon_{n}\frac
{R_{A,m}}{n},\label{deff3}\\
f_4(\overline{m},\kappa)&=&p_W(m)-p(m)-\kappa\varepsilon_{n}\frac
{R_{A,m}}{n},\label{deff4}\\
f_5(\overline{m},\kappa)&=&\delta(m,m')-\kappa\varepsilon_{n}\frac
{R_{A,m}\vee
R_{A,m'}}n.\label{deff5}
\end{eqnarray}

We will use the following fact.\vspace*{6pt}

\setcounter{fact}{-1}
\begin{fact}\label{fact0} The resampling penalty $\operatorname
{pen}_W(m,C)$ defined
in (\ref{defRespen}) satisfies
\[
\operatorname{pen}_W(m,C)=2C\tilde{C}^{-1}_{W}p_W(m).
\]
\end{fact}
\begin{pf} Let $(\psi_{\lambda})_{\lambda\in\Lambda_m}$ be an
orthonormal basis of $S_m$. Recall that $\widehat{s}_{A,m}^W=\sum
_{\lambda\in
\Lambda_m}(P_A^W\psi_{\lambda})\psi_{\lambda},$ so that
\[
\widehat{s}_{A,m}^W-\overline{W}\widehat{s}_{A,m}=\sum_{\lambda\in
\Lambda_m}(\nu^W_A\psi
_{\lambda})\psi_{\lambda}.
\]
Hence, $\nu_A^W(\widehat{s}_{A,m}^W-\overline{W}\widehat
{s}_{A,m})=\sum_{\lambda\in\Lambda
_m}(\nu_A^W\psi_{\lambda})^2.$

We conclude the proof showing that $\mathbb{E}_W(\nu_A^W(\overline
{W}\widehat{s}
_{A,m}))=0$, hence
\[
\frac{p_W(m)}{\tilde{C}_{W}}=\mathbb{E}_W\bigl(\nu_A^W(\widehat
{s}_{A,m}^W-\overline{W}\widehat{s}_{A,m})
\bigr)=\mathbb{E}_W(\nu_A^W(\widehat{s}_{A,m}^W))=\frac{\operatorname
{pen}_W(m,C)}{2C}.
\]
Since $W_0,\ldots,W_{p-1}$ are independent of $X_1,\ldots,X_n$,
\[
\mathbb{E}_W(\nu_A^W(\overline{W}\widehat{s}_{A,m}))=\frac
1{p^{2}}\sum
_{i,j=0}^{p-1}L_{q}(\psi_\lambda)(A_i)L_{q}(\psi_\lambda
)(A_j)\mathbb{E}
_W\bigl(W_i\overline{W}-(\overline{W})^2\bigr).
\]
Then, by exchangeability of the weights,
\begin{eqnarray*}
\mathbb{E}_W\bigl(W_i\overline{W}-(\overline{W})^2\bigr)&=&\frac1p\biggl(\mathbb
{E}(W_i^2)+\sum_{j\neq i}\mathbb{E}
(W_iW_j)\biggr)\\
&&{}-\frac1{p^2}\biggl(\sum_{i}\mathbb{E}(W_i^2)+\sum_{i\neq j}\mathbb{E}(W_iW_j)
\biggr)=0.
\end{eqnarray*}
\upqed\end{pf}

\subsection{\texorpdfstring{Proof of Theorem \protect\ref{theoRestau}}{Proof of Theorem 3.1}}
The proof is based on the following lemma, whose proof is given in \citet{le11}.

\begin{lemma}\label{lemconcentrationtau}
Let $X_1,\ldots,X_n$ be a strictly stationary sequence of real valued
random variables with common density $s$ and let $(S_m)_{m \in\mathcal
{M}_n}$ be a
collection of regular wavelet spaces [\textup{W}] satisfying assumptions
(\textup{H3}), (\textup{H4}). Let~$p$, $q$ be two integers satisfying $2pq=n$ and
$\frac12\sqrt{n}(\ln n)^2\leq p\leq\sqrt{n}(\ln n)^2.$

Assume that there exists $\theta>5$ such that $X_1,\ldots,X_n$ are
arithmetically [$\operatorname{AR}(\theta)$] $\tau$-mixing and satisfy S($\tau$,W).
There exist constants $\kappa_1$, $\kappa_2$, such that, for all
$i=1,\ldots,5$, for all $\overline{m}$ in $\mathcal{M}$,
%
\begin{equation}
\mathbb{E}\Bigl(\sup_{\overline{m}\in\mathcal{M}}( f_{i}(\overline
{m},\kappa_{1}) )_+\Bigr)\leq
\frac{\kappa_2}n.\label{eqconcpDtau1}
\end{equation}
\end{lemma}
It comes from Fact \ref{fact0} and the equality $2C\tilde{C}_{W}=\kappa(C)+2$
that, for all $m$ in~$\mathcal{M}_{n}$,
%
\begin{equation}
\operatorname{pen}_W(m,C)-\bigl(2+\kappa(C)\bigr)p(m)=2C\tilde
{C}_{W}^{-1}\bigl(p_{W}(m)-p(m)\bigr)\label{pen-p}.
\end{equation}
Hence, from (\ref{eqconcpDtau1}) with $i=3,4$, $\operatorname
{pen}_W(m,C)$ satisfies
conditions (\ref{condpenpastropgrande}) and~(\ref{condpenassezgrande}) of Theorem~\ref{theoSlope2tau} with $\delta
_{+}=-\delta_{-}=\kappa(C)$ and $\varepsilon=2\kappa_{1}C\tilde{C}_{W}
^{-1}\varepsilon_{n}$. Theorem~\ref{theoRestau} follows
from~(\ref{eqOISlopetau}).\vfill\eject

\subsection{\texorpdfstring{Proof of Theorem \protect\ref{theoResbeta}}{Proof of Theorem 4.1}}
The proof is based on
the following lemma whose proof is given in additional material.
\begin{lemma}\label{lemconcBeta}
Let $\theta>1$ and let $(X_n)_{n\in\mathbb{Z}}$ be an arithmetically
[$\operatorname{AR}(\theta)$] $\beta$-mixing process satisfying $\mathrm{S}(\beta
)$. Let
$(S_m)_{m \in\mathcal{M}_n}$ be a collection of linear spaces
satisfying assumptions
(\textup{H1})--(\textup{H4}). Let $p,q$ such that $2pq=n$, $\sqrt{n}(\ln
n)^2/2\leq p\leq\sqrt{n}(\ln n)^2$. There exist constants $\kappa_1$,
$\kappa_2$ and an event $\Omega_n$ satisfying
\[
\mathbb{P}(\Omega_n)\geq1-\kappa_2\biggl(\frac{(\ln n)^{2(1+\theta
)}}{n^{\theta/2}}\vee\frac1{n^2}\biggr),
\]
such that, on $\Omega_{n}$,
%
\begin{equation}
\forall\overline{m}\in\mathcal{M}, \forall i=1,\ldots,5\qquad
f_{i}(\overline{m})\leq0.
\end{equation}
\end{lemma}
Hence, from (\ref{eqconcpDtau1}) with $i=3,4$, $\operatorname
{pen}_W(m,C)$ satisfies
condition (\ref{condbonnePen}) of Theorem~\ref{theoSlope2beta} with
$\delta_{+}=-\delta_{-}=\kappa(C)$ and $\varepsilon=2\kappa
_{1}C\tilde{C}_{W}
^{-1}\varepsilon_{n}$. Theorem \ref{theoResbeta} follows
from~(\ref{eqOISlopebeta}).\vspace*{8pt}

\subsection{\texorpdfstring{Proof of Theorems \protect\ref{theoSlope1tau} and \protect\ref{theoSlope1beta}}
{Proof of Theorems 3.2 and 4.2}}
It is sufficient to prove the results for
sufficiently large $n$ since we can increase the constant $\kappa_{2}$
if necessary. Let $m_o$ be a model such that $R_{A,m_o}=R_n$. Now, by
definition, $\widehat{m}$ minimizes among $\mathcal{M}_n$ the following criterion:
\[
\mathrm{Crit}(m)=\|\widehat{s}_{A,m}\|^2-2P_A\widehat
{s}_{A,m}+\operatorname{pen}(m)+\|s\|^2+2\nu_A(s_{m_o}).
\]

\begin{fact}\label{fact1} For all $m$ in $\mathcal{M}_n$,
\[
\mathrm{Crit}(m)=\|s_m-s\|^2+\operatorname{pen}(m)-p(m)+2\nu_A(s_{m_o}-s_m).
\]
\end{fact}
\begin{pf}
 Recalling that $\|s-\widehat{s}_{A,m}\|^2=\|\widehat
{s}_{A,m}\|
^2-2P\widehat{s}
_{A,m}+\|s\|^2$ and that $(P_A-P)(\widehat{s}_{A,m}-s_m)=\|\widehat
{s}_{A,m}-s_m\|
^2=p(m)$, we have,
\begin{eqnarray*}
\mathrm{Crit}(m)&=&\|s-\widehat{s}_{A,m}\|^2-2\nu_A(\widehat
{s}_{A,m}-s_m)+2\nu
_A(s_{m_o}-s_m)+\operatorname{pen}(m)\\
&=&(\|s-\widehat{s}_{A,m}\|^2-\|\widehat{s}_{A,m}-s_m\|
^2)-p(m)+\operatorname{pen}(m)+2\nu_A(s_{m_o}-s_m).
\end{eqnarray*}
We conclude the proof with the Pythagoras equality.
\end{pf}

\begin{fact}\label{fact2} For all $m$ in $\mathcal{M}_n$, for all constants
$\kappa_1$,
\begin{eqnarray*}
(1+2\kappa_1\varepsilon_{n})\frac{2D_{A,m}}{n}&\geq&-\mathrm{Crit}(m)+(1-2\kappa_1\varepsilon_{n})\|s-s_m\|^2\\
&&{}-\sup_{m \in\mathcal{M}_n}(f_{1}(m,\kappa_{1}))-\sup_{(m,m')\in
\mathcal{M}_n^2}(f_{5}((m,m'),\kappa_{1})).
\end{eqnarray*}
\end{fact}
\begin{pf} From Fact \ref{fact1}, for all $m$ in $\mathcal{M}_n$, for all
$\kappa_1$,
since $\operatorname{pen}(m)\geq0$,
\begin{eqnarray*}
\mathrm{Crit}(m)\geq\|s_m-s\|^2-f_{1}(m,\kappa_{1})-\frac
{2D_{A,m}}n-2\kappa_1\varepsilon_{n}\frac
{R_{A,m}}n-f_{5}((m_{o},m),\kappa_{1}).
\end{eqnarray*}
We conclude the proof using that $R_{A,m}=n\|s-s_m\|^2+2D_{A,m}$.
\end{pf}

\begin{fact}\label{fact3} For all $m$ in $\mathcal{M}_n$, for all constants
$\kappa_1$,
\begin{eqnarray*}
(\delta-4\kappa_1\varepsilon_{n})\frac{D_{A,m}}{n}&\leq&-\mathrm{Crit}(m)+(1+2\kappa_1\varepsilon_{n})\|s-s_m\|^2\\
&&{}+\sup_{m \in\mathcal{M}_n}(f_{2}(m,\kappa_{1}))+\sup_{(m,m')\in
\mathcal{M}_n^2}
(f_{5}((m,m'),\kappa_{1})).
\end{eqnarray*}
\end{fact}
\begin{pf} From Fact \ref{fact1}, for all $m$ in $\mathcal{M}_n$, for all
$\kappa_1$,
since $\operatorname{pen}(m)\leq(2-\delta)D_{A,m}/n$,
\begin{eqnarray*}
\mathrm{Crit}(m)\leq\|s_m-s\|^2+f_{2}(m,\kappa_{1})-\delta\frac
{D_{A,m}}{n}+2\kappa_1\varepsilon_{n}\frac
{R_{A,m}}n+f_{5}((m,m_{o}),\kappa_{1}).
\end{eqnarray*}
We conclude the proof using that $R_{A,m}=n\|s-s_m\|^2+2D_{A,m}$.
\end{pf}

From Fact \ref{fact2}, we have, for all $\kappa_1$,
\begin{eqnarray*}
(1+2\kappa_1\varepsilon_{n})\frac{2D_{A,\widehat{m}}}{n}&\geq&
-\mathrm{Crit}(\widehat{m})+(1-2\kappa_1\varepsilon_{n})\|s-s_{\hat
{m}}\|^2\\
&&{}-\sup_{m \in\mathcal{M}_n}(f_{1}(m,\kappa_{1}))-\sup_{(m,m')\in
\mathcal{M}_n^2}
(f_{5}((m,m'),\kappa_{1})).
\end{eqnarray*}
Let us now consider a model $m^*$ such that $D_{A,m^*}=D_n^*$. By
definition of $\widehat{m}$, we have $\mathrm{Crit}(\widehat{m})\leq\mathrm{Crit}(m^*).$
Hence, from Fact \ref{fact3}, we deduce that
%
\begin{eqnarray}\label{eqfinPreuveSlope1tau}
(1+2\kappa_1\varepsilon_{n})\frac{2D_{A,\widehat{m}}}{n}&\geq&-\mathrm{Crit}(m^*)+(1-2\kappa_1\varepsilon_{n})\|s-s_{\widehat{m}}\|
^2\nonumber\\
&&{}-\sup_{m \in\mathcal{M}_n}(f_{1}(m,\kappa_{1}))-\sup_{(m,m')\in
\mathcal{M}_n^2}
(f_{5}((m,m'),\kappa_{1})).
\nonumber
\\[-8pt]
\\[-8pt]
\nonumber
&\geq&\biggl(\delta-4\kappa_1\varepsilon_{n}-(1+2\kappa_1\varepsilon
_{n})\frac{\Vert s-s_{m^*}\Vert^{2}}{D_{n}^{*}}\biggr)\frac
{D_{A,m^*}}{n}\\
&&{}+(1-2\kappa_1\varepsilon_{n})\|s-s_{\widehat{m}}\|^2-4\sup_{i\in
\{1,2,5\}, \overline{m}\in\mathcal{M}}(f_{i}(\overline{m},\kappa
_{1})).\nonumber
\end{eqnarray}
From Lemma~\ref{lemconcentrationtau}, there exist $\kappa_{1}$ and
$\kappa_{2}$ such that
\[
\mathbb{E}\Bigl( 4\sup_{i\in\{1,2,5\}, \overline{m}\in\mathcal
{M}}(f_{i}(\overline{m},\kappa _{1}))_{+} \Bigr)\leq\frac{\kappa_{2}}n.
\]
From Lemma~\ref{lemconcBeta}, there exists $\kappa_{1}$ such that, on
$\Omega_{n}$,
\[
4\sup_{i\in\{1,2,5\}, \overline{m}\in\mathcal{M}}(f_{i}(\overline
{m},\kappa_{1})
)\leq0.
\]
Now, assume that $n$ is sufficiently large to ensure that
\[
4\kappa_1\varepsilon_{n}\leq\frac{\delta}4\leq\frac14,\qquad \frac{n\|
s-s_{m^*}\|^2}{D_n^*}\leq\frac{2\delta}9.
\]
Then, taking the expectation in (\ref{eqfinPreuveSlope1tau}), we
obtain that
\[
\frac{9\mathbb{E}(D_{A,\widehat{m}})}{8n}\geq\frac{\delta}2\frac
{D_n^*}{n}-\frac{\kappa_4}{n}.
\]
Hence, (\ref{SME}) is proved for $n$ sufficiently large.

Moreover, on $\Omega_{n}$, we have
\[
\frac{9D_{A,\widehat{m}}}{8n}\geq\frac{\delta}2\frac{D_n^*}{n}.
\]
Hence, the first inequality of (\ref{eqpasdOI}) is proved for $n$
sufficiently large.
(\ref{SME2}) and the second inequality of (\ref{eqpasdOI}) follow from
the inequality
\[
\|s-\tilde{s}_A\|^2\geq(1-\kappa_1\varepsilon_{n})\frac
{R_{A,\widehat{m}}}{n}-f_{2}(\widehat{m},\kappa_{1}).
\]
From Lemma~\ref{lemconcentrationtau}, there exist constants $\kappa
_1$, $\kappa_2$, such that $\mathbb{E}(f_{2}(\widehat{m},\kappa
_{1}))\leq
\kappa_{2}/n$. We choose $n$ sufficiently large to ensure that $\kappa
_1\varepsilon_{n}\leq1/2$, we use (\ref{SME}) and we obtain that there
exists a constant $\kappa$ such that
\[
\mathbb{E}(\|s-\tilde{s}_A\|^2)\geq\frac{2\delta}9\frac
{D_n^*-\kappa}{n}.
\]
We conclude the proof of (\ref{SME2}) with the following fact.
\begin{fact}\label{fact4}
\[
\frac{R_n}n\geq\frac{16}{17}\mathbb{E}\Bigl(\inf_{m \in\mathcal
{M}_n}\|s-\widehat{s}_{A,m}\|^2
\Bigr)-\frac{\kappa}n,
\]
thus
\[
\frac{D_n^*}n\geq\frac{16D_n^*}{17 R_n}\biggl(\mathbb{E}\Bigl(\inf
_{m \in\mathcal{M}_n}\|s-\widehat{s}_{A,m}\|^2\Bigr)-\frac{\kappa}n\biggr).
\]
\end{fact}
\begin{pf} Let $\kappa_{1}$ be the constant
previously defined,
\[
\inf_{m \in\mathcal{M}_n}\|s-\widehat{s}_{A,m}\|^2\leq(1+\kappa
_{1}\varepsilon_{n})\inf_{m \in\mathcal{M}_n
}\biggl\{ \frac{R_{A,m}}n \biggr\}+\sup_{m \in\mathcal{M}_n}f_{1}(m,\kappa_{1}).
\]
We conclude the proof with Lemma~\ref{lemconcentrationtau}.
\end{pf}

We use the first inequality of (\ref{eqpasdOI}) and we obtain that, on
$\Omega_{n}$,
\[
\|s-\tilde{s}_A\|^2\geq\frac{2\delta}9\frac{D_n^*}{n}.
\]
We conclude the proof of Theorem~\ref{theoSlope1beta}, saying that,
on $\Omega_n$, we have
\begin{eqnarray*}
\frac{R_n}n&=&\inf_{m\in\mathcal{M}_n}\biggl\{ \|s-s_m\|^2+\frac
{2D_{A,m}}n \biggr\}\geq
(1-\kappa_1\varepsilon_n)\inf_{m\in\mathcal{M}_n}\{ \|s-s_m\|
^2+p(m) \}\\
&=&(1-\kappa_1\varepsilon_n)\inf_{m\in\mathcal{M}_n}\|s-\widehat
{s}_{A,m}\|
^2\geq\frac{15}{16}\inf_{m\in\mathcal{M}_n}\|s-\widehat{s}_{A,m}\|^2.
\end{eqnarray*}
Thus,
\[
\|\tilde{s}_A-s\|^2\geq\frac{2\delta}{9}\frac{D_n^*}{R_n}\frac
{R_n}n\geq\frac{\delta}{9}\frac{D_n^*}{R_n}\inf_{m\in\mathcal
{M}_n}\|s-\widehat{s}
_{A,m}\|^2.\vspace*{8pt}
\]

%
%
\subsection{\texorpdfstring{Proofs of Theorems \protect\ref{theoSlope2tau} and \protect\ref{theoSlope2beta}}
{Proofs of Theorems 3.3. and 4.3}}
As in the previous proof, it is sufficient to obtain the results for
sufficiently large $n$. Let us first prove the oracle inequalities. Let
$\kappa_{1}$ be a constant to be chosen later. Let $\Omega_{n}$ be the
set defined on Lemma~\ref{lemconcBeta}. The key point to prove oracle
inequalities is the following fact.
\begin{fact}\label{fact5} For all $m$ in $\mathcal{M}_n$, for all real numbers
$\delta_{-}$,
$\delta_{+}$ and for all nonnegative reals $x,y$,
%
\begin{eqnarray}\label{eqRemTerm1}
\qquad&&\bigl[\bigl(1\wedge(1-\delta_{-})\bigr)-x-y\bigr]\|s-\tilde{s}_A\|^2\nonumber\\
&&\qquad\leq\bigl[\bigl(1\vee
(1+\delta
_{+})\bigr)+x+y\bigr]\|s-\widehat{s}_{A,m}\|^2
\\
&&\qquad\quad{}\!+\sup_{m\in\mathcal{M}_{n}}\{ \operatorname{pen}(m)-(2+\delta
_{+})\|\widehat{s}_{A,m}-s_m\|^2-x\| s-\widehat{s}_{A,m}\|^2 \}
_{+}\nonumber\\
&&\qquad\quad{}\!+\sup_{m\in\mathcal{M}_{n}}\{ (2-\delta_{-}) \|\widehat
{s}_{A,m}-s_m\|^2-\operatorname{pen}(m)-x\| s-\widehat{s}_{A,m}\|^2 \}
_{+}\label{eqRemTerm2}\\
&&\qquad\quad{}\!+2\!\sup_{(m,m')\in\mathcal{M}_{n}^{2}}\{ \nu_A(s_{m'}-s_m)-y(\|
s-\widehat{s}_{A,m}\| ^2+\|s-\widehat{s}_{A,m'}\|^2) \}_{+}.\label{eqRemTerm3}
\end{eqnarray}
\end{fact}
\begin{pf} By definition of $\tilde{s}_A$, for all $m$ in
$\mathcal{M}_n$,
we have
\[
\|\tilde{s}_A\|^2-2P_A\tilde{s}_A+\operatorname{pen}(\widehat{m})+\|s\|
^2\leq\|\widehat{s}
_{A,m}\|^2-2P_A\widehat{s}_{A,m}+\operatorname{pen}(m)+\|s\|^2.
\]
Now, for all $m$ in $\mathcal{M}_n$, since $\|\widehat{s}_{A,m}-s\|
^2=\|\widehat{s}_{A,m}\|
^2-2P\widehat{s}_{A,m}+\|s\|^2$,
\[
\|\widehat{s}_{A,m}\|^2-2P_A\widehat{s}_{A,m}+\|s\|^2=\|\widehat
{s}_{A,m}-s\|^2-2(P_A-P)\widehat{s}_{A,m}.
\]
Thus, for all $m$ in $\mathcal{M}_n$,
\[
\|\tilde{s}_A-s\|^2-2(P_A-P)\tilde{s}_A+\operatorname{pen}(\hat
{m})\leq\|\widehat{s}
_{A,m}-s\|^2-2(P_A-P)\widehat{s}_{A,m}+\operatorname{pen}(m).
\]
For all $m$ in $\mathcal{M}_n$, since $(P_A-P)(\widehat
{s}_{A,m}-s_m)=\|\widehat{s}_{A,m}-s\|^2$,
\[
2(P_A-P)\widehat{s}_{A,m}=2\|s_m-\widehat{s}_{A,m}\|^2+2(P_A-P)s_m.
\]
This yields
\begin{eqnarray*}
\|s-\tilde{s}_A\|^2&\leq&\|s-\widehat{s}_{A,m}\|^2+\operatorname
{pen}(m)-2\|\widehat{s}_{A,m}-s_m\|
^2\\
&&{}+2\|\widehat{s}_{A,\widehat{m}}-s_{\widehat{m}}\|^2-\operatorname
{pen}(\widehat{m})+2\nu_A(s_{\widehat{m}}-s_m).
\end{eqnarray*}
We add $-[(\delta_{-}\vee0)+(x+y)]\|\tilde{s}_A-s\|^2$ to the
left-hand side of the previous inequality and $-\delta_{-}\|\tilde
{s}_A-s_{\widehat{m}}\|^2-(x+y)\|\tilde{s}_A-s\|^2+[(\delta_{+}\vee
0)+x+y]\|
s-\widehat{s}_{A,m}\|^2-\delta_{+}\|\widehat{s}_{A,m}-s_m\|^2-(x+y)\|
s-\widehat{s}_{A,m}\|
^2$ to the right-hand side. This is valid because, for all $m$ in
$\mathcal{M}
_n$, for all reals $\delta$,
\[
[(\delta\vee0)+x+y]\|\widehat{s}_{A,m}-s\|^2\geq\delta\|\widehat
{s}_{A,m}-s_{m}\|
^2+(x+y)\|\widehat{s}_{A,m}-s\|^2.
\]
We obtain
\begin{eqnarray*}
&&\bigl[\bigl(1\wedge(1-\delta_{-})\bigr)-x-y\bigr]\|s-\tilde{s}_A\|^2\\
&&\qquad\leq\bigl[\bigl(1\vee
(1+\delta
_{+})\bigr)+x+y\bigr]\|s-\widehat{s}_{A,m}\|^2\\
&&\qquad\quad{}+\operatorname{pen}(m)-(2+\delta_{+})\|\widehat{s}_{A,m}-s_m\|^2-x\|
\widehat{s}_{A,m}-s\|^2\\
&&\qquad\quad{}+(2-\delta_-)\|\widehat{s}_{A,\widehat{m}}-s_{\widehat{m}}\|
^2-\operatorname{pen}(\widehat{m})-x\|\widehat{s}_{A,\widehat
{m}}-s\|
^2\\
&&\qquad\quad{}+2\nu_A(s_{\widehat{m}}-s_m)-y\|\widehat{s}_{A,m}-s\|^2-x\|\widehat
{s}_{A,\widehat{m}}-s\|
^2.
\end{eqnarray*}
\upqed\end{pf}
We will also use the following fact.

\begin{fact}\label{fact6} For all reals $\kappa$ such that $\kappa
\varepsilon
_{n}\leq1/2$,
\[
\frac{R_{A,m}}n\leq2\|s-\widehat{s}_{A,m}\|^2+2\{ f_{2}(m,\kappa) \}_+.
\]
\end{fact}
\begin{pf} We write
\[
\frac{R_{A,m}}n=\frac1{1-\kappa\varepsilon_{n}}\biggl( \frac
{R_{A,m}}n-\|s-\widehat{s}_{A,m}\|^2-\kappa\varepsilon_{n}\frac
{R_{A,m}}n \biggr)+\frac1{1-\kappa\varepsilon_{n}}\|s-\widehat{s}_{A,m}\|^2.
\]
We use that $\kappa\varepsilon_{n}\leq1/2$ and that
$R_{A,m}=2D_{A,m}+n{s-s_m}^2$ to conclude the proof.
\end{pf}

\textit{Control of \textup{(\ref{eqRemTerm1})}.}
Assume that $n$ is sufficiently large to ensure that $\kappa
_{1}\varepsilon_{n}\leq1/2$. We have, from Fact \ref{fact6},
\begin{eqnarray*}
&&\operatorname{pen}(m)-(2+\delta_{+})p(m)-2\varepsilon{\Vert\widehat
{s}_{A,m}-s\Vert }^{2}\\
&&\qquad\leq\operatorname{pen}(m)-(2+\delta_{+})p(m)-\varepsilon\frac
{R_{A,m}}n+2\varepsilon
\{ f_{2}(m,\kappa_{1}) \}_+.
\end{eqnarray*}
Applying Lemma~\ref{lemconcentrationtau}, we obtain constants
$\kappa
_{1}$ and $\kappa_{2}$ such that
\[
\mathbb{E}\Bigl(\sup_{m \in\mathcal{M}_n} \{ f_{2}(m,\kappa_{1}) \}
_{+}\Bigr)\leq\frac{\kappa_{2}}n.
\]
Applying Lemma~\ref{lemconcBeta}, we obtain a constant $\kappa_{1}$
such that, on $\Omega_{n}$,
\[
\sup_{m \in\mathcal{M}_n} \{ f_{2}(m,\kappa_{1}) \}_{+}\leq0.
\]
Moreover, (\ref{condpenpastropgrande}) ensures that
\[
\mathbb{E}\biggl(\sup_{m \in\mathcal{M}_n} \biggl\{ \operatorname
{pen}(m)-(2+\delta_{+})p(m)-\varepsilon\frac {R_{A,m}}n \biggr\}_{+}\biggr)\leq
\frac{\kappa}n.
\]
On $\Omega_{\operatorname{pen}}$, we have
\[
\sup_{m \in\mathcal{M}_n} \biggl\{ \operatorname{pen}(m)-(2+\delta
_{+})p(m)-\varepsilon\frac {R_{A,m}}n \biggr\}_{+}\leq0.
\]
We choose $x=2\varepsilon$. We obtain that, for Theorem~\ref{theoRestau}, the expectation of (\ref{eqRemTerm1}) is upper bounded
by $\kappa n^{-1}$ and for Theorem~\ref{theoResbeta}, the term (\ref{eqRemTerm1}) is equal to~$0$ on $\Omega_{n}\cap\Omega
_{\operatorname{pen}}$.

\textit{Control of \textup{(\ref{eqRemTerm2})}.} Assume that $n$ is sufficiently
large to ensure that $\kappa_{1}\varepsilon_{n}<1/2$, we deduce from
Fact \ref{fact6} that
\begin{eqnarray*}
&&(2-\delta_{-})p(m)-\operatorname{pen}(m)-2\varepsilon{\Vert\widehat
{s}_{A,m}-s\Vert}^{2}\\
&&\qquad\leq(2-\delta_{-})p(m)-\operatorname{pen}(m)-\varepsilon\frac
{R_{A,m}}n+2\varepsilon
\{ f_{2}(m,\kappa_{1}) \}_+.
\end{eqnarray*}
Applying Lemma~\ref{lemconcentrationtau}, we obtain constants
$\kappa
_{1}$ and $\kappa_{2}$ such that
\[
\mathbb{E}\Bigl(\sup_{m \in\mathcal{M}_n} \{ f_{2}(m,\kappa_{1}) \}
_{+}\Bigr)\leq\frac{\kappa_{2}}n.
\]
Applying Lemma~\ref{lemconcBeta}, we obtain a constant $\kappa_{1}$
such that, on $\Omega_{n}$,
\[
\sup_{m \in\mathcal{M}_n} \{ f_{2}(m,\kappa_{1}) \}_{+}\leq0.
\]
Moreover, (\ref{condpenassezgrande}) ensures that
\[
\mathbb{E}\biggl(\sup_{m \in\mathcal{M}_n} \biggl\{ (2-\delta
_{-})p(m)-\operatorname{pen}(m)-\varepsilon\frac {R_{A,m}}n \biggr\}
_{+}\biggr)\leq\frac{\kappa}n.
\]
On $\Omega_{\operatorname{pen}}$, we have
\[
\sup_{m \in\mathcal{M}_n} \biggl\{ (2-\delta_{-})p(m)-\operatorname
{pen}(m)-\varepsilon\frac {R_{A,m}}n \biggr\}_{+}\leq0.
\]
We choose $x=2\varepsilon$. We obtain that, for Theorem~\ref{theoRestau}, the expectation of (\ref{eqRemTerm2}) is upper bounded
by $\kappa n^{-1}$ and for Theorem~\ref{theoResbeta}, the term (\ref{eqRemTerm2}) is equal to~$0$ on $\Omega_{n}\cap\Omega
_{\operatorname{pen}}$.

\textit{Control of \textup{(\ref{eqRemTerm3})}.} Let $m,m'$ in $\mathcal
{M}_n$ and let
$m_{s}$ be the index such that $R_{A,m_{s}}=R_{A,m}\vee R_{A,m'}$ and
let $\kappa_{1}$ be a constant to be chosen later. Assume that $n$ is
sufficiently large to ensure that $\kappa_{1}\varepsilon_{n}\leq1/2$.
It comes from Fact \ref{fact6} that
\begin{eqnarray*}
\delta(m,m')&=&f_{5}((m,m'),\kappa_{1})+\kappa_{1}\varepsilon
_{n}\frac
{R_{A,m_{s}}}n\\
&\leq& f_{5}((m,m'),\kappa_{1})+2\kappa_{1}\varepsilon_{n}{\Vert\widehat{s}
_{A,m_{s}}-s\Vert}^{2}+2\kappa_{1}\varepsilon_{n}\{ f_{2}(m_{s},\kappa
_{1}) \}_+.
\end{eqnarray*}
We deduce from Lemma~\ref{lemconcentrationtau} that there exist
$\kappa_{1}$ and $\kappa_{2}$ such that
\begin{eqnarray*}
&&\mathbb{E}\Bigl( \sup_{(m,m')\in\mathcal{M}_{n}^{2}}\{ \delta
(m,m')-2\kappa _1\varepsilon_{n}(\|\widehat{s}_{A,m}-s\|^2+\|
\widehat{s}_{A,m'}-s\|^2) \} \Bigr)\\
&&\qquad\leq\mathbb{E}\Bigl( \sup_{(m,m')\in\mathcal{M}_{n}^{2}}\{
f_{5}((m,m'), \kappa_{1})+2\kappa _{1}\varepsilon_{n}f_{2}(m_{s},\kappa_{1}) \}
_{+} \Bigr)\leq\frac{\kappa_{2}}n.
\end{eqnarray*}
Applying Lemma~\ref{lemconcBeta}, we obtain a constant $\kappa_{1}$
such that, on $\Omega_{n}$,
\[
\sup_{(m,m')\in\mathcal{M}_{n}^{2}}\{ \delta(m,m')-2\kappa
_1\varepsilon_{n}(\| \widehat{s}_{A,m}-s\|^2+\|\widehat{s}_{A,m'}-s\|
^2) \}\leq0.
\]
\textit{Conclusion of the proofs.} We use Fact \ref{fact5} with $x=2\varepsilon$
and $y=2\kappa_1\varepsilon_{n}$. We take the expectation for the proof
of Theorem~\ref{theoRestau}, we have obtained that the expectation of
the remainder terms (\ref{eqRemTerm1})--(\ref{eqRemTerm3}) are upper
bounded by $\kappa n^{-1}$ for a sufficiently large $n$. For the proof
of Theorem~\ref{theoResbeta}, we have obtained that the remainder
terms~(\ref{eqRemTerm1})--(\ref{eqRemTerm3}) with $x=2\varepsilon$
and $y=2\kappa_1\varepsilon_{n}$ are equal to $0$ on $\Omega_{n}\cap
\Omega_{\operatorname{pen}}$ when~$n$ is sufficiently large. As
explained in the
beginning of the proof, this is sufficient to conclude the proof of
(\ref{eqOISlopetau}) and~(\ref{eqOISlopebeta}).

Let us prove (\ref{eqDhmpetit}). Let $\kappa_1<1/(2\varepsilon_{n})$,
from Fact \ref{fact6} and (\ref{eqOISlopetau}), we have
\begin{eqnarray*}
\frac{K_n}n\mathbb{E}(2D_{A,\widehat{m}})&\leq& K_n\biggl( \mathbb
{E}(p(\widehat{m}))+\mathbb{E} (f_{2}(\widehat{m},\kappa
_{1}))+\kappa_{1}\varepsilon_{n}\mathbb{E}\biggl( \frac {R_{A,\widehat
{m}}}n \biggr) \biggr)\\
&\leq&(1+2\kappa_{1}\varepsilon_{n})\mathbb{E}(( f_{2}(\widehat
{m},\kappa _{1}) )_+)+(1+2\kappa_{1}\varepsilon_{n})K_n\mathbb{E}(\|
s-\tilde{s}_A\|
^2)\\
&\leq&2\mathbb{E}(( f_{2}(\widehat{m},\kappa_{1}) )_+)+2K_n\mathbb
{E}(\|s-\tilde
{s}_A\|^2)\\
&\leq&2\biggl( \mathbb{E}(( f_{2}(\widehat{m},\kappa_{1}) )_+)+\frac
{\kappa}n \biggr)+2R_n.
\end{eqnarray*}
We used that, by definition $K_n\leq1$. We conclude the proof with
Lemma~\ref{lemconcentrationtau}.

In order to get the bound on $D_{A,\widehat{m}}$ in Theorem \ref{theoSlope2beta}, we use that, on $\Omega_n\cap\Omega
_{\operatorname{pen}}$, (\ref{eqOISlopebeta}) holds and there exists a constant $\kappa_1$ such
that, $\kappa_1\varepsilon_n<1/2$ satisfying
\begin{eqnarray*}
K_n\frac{2D_{A,\widehat{m}}}n&\leq&\frac{K_n}{1-\kappa_1\varepsilon_n}
\bigl(\|s-s_{\widehat{m}}\|^2+p^*(\widehat{m})\bigr)=\frac{K_n}{1-\kappa
_1\varepsilon_n}\|s-\tilde{s}_A\|^2\\
&\leq&\frac{1}{1-\kappa_1\varepsilon_n}\inf_{m\in\mathcal{M}_n}\|
s-\widehat{s}_{A,m}\|
^2\leq\frac{1+\kappa_1\varepsilon_n}{1-\kappa_1\varepsilon
_n}\frac
{R_n}n\leq3\frac{R_n}n.
\end{eqnarray*}
%
%
\begin{appendix}
\section*{Appendix}\vspace*{6pt}
\label{SectionAppendix}
We present in this section some classical collections of models and
prove that they satisfy (H4).

\textit{Regular histograms:} Let $d$ be an integer and let $S_{d}$ be the
space of functions~$t$ constant on all the intervals
$([k/d,(k+1)/d))_{k\in\mathbb{Z}}$. $S_{d}$ is called the space of regular
histograms with size $1/d$. The family $(\psi_{k})_{k\in\mathbb
{Z}}$, where,
for all~$k$ in $\mathbb{Z}$, $\psi_{k}=\sqrt{d}\mathbf
{1}_{[k/d,(k+1)/d)}$ is an
orthonormal basis of $S_{d}$. Let $B_{d}=\{ t\in S_{d}, {t}^{2}\leq1 \}
$. From the Cauchy--Schwarz inequality, we have
\[
\sup_{t\in B_{d}}t^{2}=\sum_{k\in\mathbb{Z}}\psi_{k}^{2}=d\mathbf
{1}_{\mathbb{R}}.
\]
Hence,
\[
b_{m}^{2}=\Bigl\Vert\sup_{t\in B_{d}}t^{2}\Bigr\Vert_{\infty}=d, \qquad P\Bigl( \sup _{t\in
B_{d}}t^{2} \Bigr)=d P( \mathbf{1}_{\mathbb{R}} )=d.
\]
(H4) holds on all the spaces $S_{d}$ with $c_{D}=1$, therefore, it
holds on the collection $(S_{d})_{d=1,\ldots,n}$ called the regular
histograms collection.

\textit{Fourier spaces:} Let $k\geq1$ be an integer and let, for all $x$
in $[0,1]$,
\[
\psi_{1,k}(x)=\sqrt{2}\cos(2\pi k x),\qquad \psi_{2,k}(x)=\sqrt{2}\sin
(2\pi
k x),\qquad \psi_{0}=\mathbf{1}_{[0,1]}.
\]
Let $\mathcal{M}_{n}=\{ 1,\ldots,n \}$ and $\forall m\in\mathcal
{M}_{n}$, let $\Lambda
_{m}=\{ 0,(1,k),(2,k),k=1,\ldots,m \}$. The space $S_{m}$, spanned by the
family $(\psi_{\lambda})_{\lambda\in\Lambda_m}$ is called the
Fourier space with harmonic
smaller than $m$ and the collection $(S_{m}, m\in\mathcal{M}_{n})$ is called
the collection of Fourier spaces. Let $B_{m}=\{ t\in S_{m}, {t}^{2}\leq
1 \}$. From the Cauchy--Schwarz inequality, for all $x$ in $[0,1]$,
\[
\sup_{t\in B_{m}}t^{2}(x)=\sum_{\lambda\in\Lambda_m}\psi_{\lambda
}^{2}(x)=1+2\sum
_{k=1}^{m}\bigl(\cos^{2}(2\pi k x)+\sin^{2}(2\pi k x)\bigr)=1+2m.
\]
Hence, if $P$ is supported in $[0,1]$,
\[
b_{m}^{2}=\Bigl\Vert\sup_{t\in B_{m}}t^{2}\Bigr\Vert_{\infty}=1+2m, P\Bigl( \sup _{t\in
B_{m}}t^{2} \Bigr)=1+2m.
\]
(H4) holds with $c_{D}=1$ on the collection of Fourier spaces when
$P$ is supported on $[0,1]$.

\textit{Wavelet spaces:} Assume that $(S_{m},m\in\mathcal{M}_{n})$
is a collection
of wavelet spaces [W]. Assume moreover that the scaling function
$\phi$ and the mother wavelet $\psi$ satisfy the following relation.
There exists a constant $K_{o}>0$ such that, for all $x$ in~$\mathbb{R}$,
\[
\frac{1}{K_{o}}\leq\sum_{k\in\mathbb{Z}}\phi^{2}(x-k)\leq K_{o},\qquad
\frac
{1}{K_{o}}\leq\sum_{k\in\mathbb{Z}}\psi^{2}(x-k)\leq K_{o}.
\]
This condition is satisfied by the Haar basis, where $\phi=\mathbf{1}
_{[0,1)}$, $\psi=\mathbf{1}_{[0,1/2)}-\mathbf{1}_{[1/2,1)}$, with $K_{o}=1$.
Then, for all $j\geq0$, we have
\[
\frac{1}{K_{o}}\leq\sum_{k\in\mathbb{Z}}\phi^{2}(2^{j}x-k)\leq
K_{o},\qquad \frac
{1}{K_{o}}\leq\sum_{k\in\mathbb{Z}}\psi^{2}(2^{j}x-k)\leq K_{o}.
\]
Let $B_{m}=\{ t\in S_{m}, {t}^{2}\leq1 \}$. From the Cauchy--Schwarz
inequality, we have
\[
\Psi_{m}(x)=\sup_{t\in B_{m}}t^{2}(x)=\sum_{\lambda\in\Lambda
_m}\psi_{\lambda}^{2}(x)=\sum
_{k\in\mathbb{Z}}2\phi^{2}(2x-k)+\sum_{j=1}^{J_{m}}2^{j}\sum_{k\in
\mathbb{Z}}\psi
^{2}(2^{j}x-k).
\]
We deduce that
\[
\frac{2^{J_{m}}}{K_{o}}\leq\frac{1}{K_{o}}\Biggl( 2+\sum
_{j=1}^{J_{m}}2^{j} \Biggr)\leq\Psi_{m}(x)\leq K_{o}\Biggl( 2+\sum
_{j=1}^{J_{m}}2^{j} \Biggr)\leq2K_{o}2^{J_{m}}.
\]
Hence, $b_{m}^{2}={\Psi_{m}}_{\infty}\leq2K_{o}2^{J_{m}}, P( \Psi
_{m} )\geq2^{J_m}/K_o.$

(H4) holds wit $c_{D}=1/(2K_{o}^{2})$ on the collection [W].
\end{appendix}

\section*{Acknowledgments} The author is very grateful to B\'{e}atrice
Laurent and Cl\'{e}mentine Prieur for their precious advice. He would
like also to thank the reviewers and Associate Editors for their
careful reading of the manuscript and helpful comments which led to an
improved presentation of the paper.

%

\begin{supplement}[id=suppA]
\stitle{Proofs of Lemmas \ref{lemconcentrationtau} and \ref{lemconcBeta}}
\slink[doi]{10.1214/11-AOS888SUPP}
\sdatatype{.pdf}
\sfilename{aos888\_suppl.pdf}
\sdescription{In the Supplementary Material, we give complete proofs
of the concentrations Lemmas \ref{lemconcentrationtau} and \ref{lemconcBeta}. We use coupling results, respectively, of \citet{Be79}
and \citet{DP05}, to build sequences of independent random variables
$(A_{0}^{*},\ldots,A_{p-1}^{*})$ approximating the sequence of blocks
$(A_{0},\ldots,A_{p-1})$, respectively in the $\beta$ and $\tau$ mixing
case. We prove concentration lemmas equivalent to Lemmas \ref{lemconcentrationtau} and \ref{lemconcBeta} for these approximating
random variables. The main tools here are the concentration
inequalities of \citet{Bo02} and \citet{KR05} for the maximum of the
empirical process. We prove finally some covariance inequalities to
evaluate the expectation of $p(m)$ and deduce the rates $\varepsilon
_{n}=(\ln n)^{-1/2}$. }
\end{supplement}

%

\printaddresses

\end{document}